\documentclass[12pt,letterpaper]{amsart}
\usepackage{amsmath}
\usepackage{graphicx}
\usepackage{latexsym,amsmath,amsthm,amssymb,amsfonts,mathrsfs}
\usepackage[utf8]{inputenc}
\usepackage{xcolor}
\usepackage{lipsum}%
\usepackage{tikz}
\usepackage{xcolor}
\usepackage{comment}
\usepackage{newtx}
\usepackage{relsize}

\usepackage{geometry}
\usepackage{enumitem}
\usepackage{musicography}
\usepackage{leftindex}
\usepackage{appendix}
\usepackage[numbers]{natbib}
\usepackage{cite}
\usepackage{url}
\usepackage{amsthm}
\usepackage{commath}
\usepackage{comment}
\providecommand{\R}{\mathbb{R}} %% Real numbers
\providecommand{\Z}{\mathbb{Z}} %% Integer numbers
\providecommand{\C}{\mathbb{C}} %% Complex numbers
\newcommand\numberthis{\addtocounter{equation}{1}\tag{\theequation}}    
\providecommand{\rn}[1]{ \textup{\uppercase\expandafter{\romannumeral#1}}} %% Roman numerals

\providecommand{\Div}{\text{div}} %% divergence

\newtheorem{theorem}{Theorem}[section]
\newtheorem{lemma}[theorem]{Lemma}

\theoremstyle{definition}
\newtheorem{definition}[theorem]{Definition}

\newtheorem{cor}[theorem]{Corollary}

\theoremstyle{remark}
\newtheorem{remark}[theorem]{Remark}

\numberwithin{equation}{section}

\title[Hamiltonian stationarity of twisted Lagrangian tori]{On Hamiltonian stationarity of twisted Lagrangian tori in ${\mathbb{C}}^2$}

\author{Jingyi Chen and Patrik Coulibaly}
\address{Department of Mathematics, 
The University of British Columbia, Vancouver, BC, Canada}
\email{jychen@math.ubc.ca}
\address{Department of Mathematics, The University of British Columbia, Vancouver, BC, Canada}
\email{cpatrik@math.ubc.ca}

\thanks{Chen is partially supported by NSERC Discovery Grant, No. 22R80062}
\date{}

\begin{document}
\maketitle

\begin{abstract}
Chekanov's exotic tori have been playing an important role in symplectic geometry as they are the only known examples of Lagrangian tori in $\C^2$ that are not Hamiltonian isotopic to a product torus. In this paper, we explore the differential geometry of a wider range of tori constructed by twisting simple closed planar curves, which include both certain product tori and Chekanov's exotic tori. In particular, we investigate the minimality of area of such twisted tori under Hamiltonian deformations and show that the only minimal twisted tori are the product ones. This tells us that Chekanov's exotic tori are not area minimal in their Hamiltonian isotopy classes. 
\end{abstract} 

\section{Introduction}\label{section:introduction}

A  Lagrangian torus in $\C^2$ that cannot be deformed into a product torus through a Hamiltonian isotopy is called \textit{exotic}. An exotic torus was first discovered by Chekanov in \citep{chekanov1996lagrangian}, cf. also Eliashberg-Polterovich \citep{eliashberg1997problem}. Chekanov's exotic tori can be constructed by lifting a certain type of simple closed plane curves to $\C^2$ with twisting, see \citep{brendel2021introduction}, \citep{chekanov2010notes} or the broader survey \citep{polterovich2024lagrangian}. 

For a simple plane curve $\gamma:S^1\to \C$ that does not pass through the origin, we define
\[
    L_{\gamma} : = \left\{\frac{1}{\sqrt{2}}\left(\gamma e^{i\alpha}, \gamma e^{-i\alpha}\right) \bigg| \hspace{5pt} \alpha \in [0,2\pi) \right\} 
 \]
 which is an immersed Lagrangian torus in $\C^2$. It is necessary to exclude curves that pass through the origin, since for such curves the construction results in a "pinched torus" that is singular at the origin. Since $L_\gamma$ is obtained by lifting $\gamma$ to the diagonal of $\C^2$ and then twisting it \citep{chekanov2010notes}, we may refer to $L_\gamma$ as a twisted torus. When $\gamma$ is contained in an open half-plane, then $L_\gamma$ is an instance of Chekanov's exotic tori in the form they are presented in \citep{eliashberg1997problem}. Therefore, when $\gamma$ is contained in an open half-plane we may refer to $L_\gamma$ as a Chekanov torus. 

Another way to construct $L_\gamma$ is by lifting a curve from a certain reduced symplectic space that can be identified with $\C^*  := \C \setminus \{0\}$. We can always find a curve $\tilde{\gamma} \subset \C^*$ that lifts to $L_\gamma$, however, $\tilde{\gamma}$ depends on the choice of the identification and it may or may not be same curve as $\gamma$. In Appendix \ref{appendix:construction}, we discuss how to find such a $\tilde{\gamma}$. Any Hamiltonian isotopy that is supported in $\C^*$ can be lifted to a Hamiltonian isotopy in $\C^2$, see Section 2 and 3 in \citep{brendel2021introduction} for more details. Also, any simple closed curve $\tilde{\gamma}$ in $\C^*$ that encloses the origin can be deformed into a circle centred at the origin through a Hamiltonian isotopy supported in $\C^*$. All such circles lift to product tori, hence each $L_\gamma$ constructed from such a $\tilde{\gamma}$ is Hamiltonian isotopic to a product torus. Similarly, all simple closed curves $\tilde{\gamma}$ in $\C^*$ that do not go around the origin and enclose the same area can be deformed into each other through a Hamiltonian isotopy supported in $\C^*$. So, after fixing the area enclosed by $\tilde{\gamma}$, the $L_\gamma$ constructed from such curves represent at most one new Hamiltonian class. No isotopy between a curve that encloses the origin and a curve that does not enclose the origin can avoid the origin. Therefore, such an isotopy is not supported in $\C^*$ and cannot be lifted to $\C^2$. This suggests that an $L_\gamma$ that corresponds to a $\tilde{\gamma}$ not enclosing the origin might be exotic. Indeed, this is the case as demonstrated by Chekanov in \citep{chekanov1996lagrangian}.

%passing through zero ----> Whitney sphere

The exotic tori reveal important symplectic geometric properties. However, their differential geometric structures seem less explored. 
In this paper, we are interested in the minimality of area of $L_\gamma$ under Hamiltonian variations, i.e. variations generated by vector fields of the form $V = J\nabla h$ for some $h\in C^\infty(L_\gamma)$ where $J$ is the standard complex structure on  $\C^2$. We say that $L_\gamma$ is Hamiltonian stationary if it is a critical point of the volume functional under Hamiltonian variations. In \citep{oh}, Y.-G. Oh investigated this variational problem and found that a Lagrangian submanifold is Hamiltonian stationary if and only if
 \begin{align}\label{eqn:E-L_general}
     \Div_g(JH) = 0,
 \end{align}
 where $H$ is the mean curvature vector of the Lagrangian submanifold and $g$ is the induced metric on it. Our main observation is as follows.

\begin{theorem}\label{theorem:main}
    $L_\gamma \subset \C^2$ is Hamiltonian stationary if and only if $\gamma  \subset \C$ is a circle centred at the origin. 
\end{theorem}

Since Chekanov's exotic tori are of the form $L_\gamma$ for a simple loop $\gamma$ that lies in an open half-plane of $\C^2$, Theorem \ref{theorem:main} yields
\begin{cor}
    Chekanov's exotic tori are not Hamiltonian stationary.
\end{cor}

It is interesting to determine whether there exist Hamiltonian stationary tori in the Hamiltonian isotopy class of a 
twisted torus and whether a product torus is the only Hamiltonian stationary one in its class. The latter relates to an open question of Y.-G. Oh on whether the product torus in $\C^2$ is area minimizing in its Hamiltonian isotopy class. 
Theorem \ref{theorem:main} says that any competitor cannot arise as $L_\gamma$. 

The outline of the paper is as follows. In Section \ref{section:DiffGeo}, we use an explicit parametrization of $L_\gamma$ to calculate some basic differential geometric quantities. 
%In Section \ref{section:Construction}, we give a brief description of the construction of the exotic tori.
In Section \ref{section:minimality}, we explore properties of Hamiltonian stationary $L_\gamma$. In Section \ref{section:Proof}, we give a proof of Theorem \ref{theorem:main}. In Appendixes, we provide a detailed account, in differential-geometric terms, on the construction and properties of the twisted tori via symplectic reduction mentioned earlier. 

\section{Differential Geometry of the twisted torus $L_\gamma$}\label{section:DiffGeo}
 Let $F:[0,2\pi)\times [0,2\pi)\to\C^2$ be given by
\[
F(\alpha,\beta)=\frac{\rho(\beta)}{\sqrt{2}} \left(e^{i\left(f(\beta)+ \alpha\right)},e^{i\left(f(\beta)-\alpha\right)}\right).
\]
Then $L_\gamma$ is parametrized by $F$ so the vectors
\begin{align*}
    e_1 &= \partial_\alpha F=  \frac{i\rho}{\sqrt{2}} \left(e^{i\left(f+ \alpha\right)},-e^{i\left(f-\alpha\right)}\right),\\
     e_2 &= \partial_\beta F=  \frac{\dot{\rho} + i \rho \dot{f}}{\sqrt{2}} \left(e^{i\left(f+ \alpha\right)},e^{i\left(f-\alpha\right)}\right) = \tau F
\end{align*}
are tangent to $L_\gamma$ where $\tau = \frac{\dot{\rho}}{\rho} + i\dot{f}$. Let us define $v:= \frac{d}{dr}\ln{\rho} = \frac{\dot{\rho}}{\rho}$ and $w := \dot{f}$. Then $\tau = v + iw$ and the induced metric $g$ on $L_\gamma$ is given by
\begin{align*}
    g_{\alpha \alpha} &= \abs{\partial_{\alpha}F}^2= \frac{\rho^2}{2}\left(\abs{e^{i(f(\beta)+\alpha +\frac{\pi}{2})}}^2 + \abs{e^{i(f(\beta) - \alpha - \frac{\pi}{2})}}^2\right)=\rho^2 \\
    g_{\alpha \beta} &= g_{\beta \alpha} = \text{Re}\left(\langle \partial_{\beta}F,\overline{\partial_{\alpha}F} \rangle \right) \\ &=\text{Re}\bigg( \frac{-i\rho^2 \tau}{2}\left\langle \left(e^{i\left(f+ \alpha\right)},e^{i\left(f-\alpha\right)}\right),\left(e^{-i\left(f+ \alpha\right)},-e^{-i\left(f-\alpha\right)}\right)  \right\rangle \bigg)\\
    &=\text{Re}\bigg( \frac{-i\rho^2 \tau}{2}(1 -1) \bigg) = 0\\
    g_{\beta \beta} &= \abs{\partial_{\beta}F}^2= \abs{\tau}^2\abs{F}^2 = (v^2+w^2)\rho^2 = \dot{\rho}^2 + \rho^2 \dot{f}^2,
    \end{align*}
    i.e.
    \begin{align*}
         g= \rho^2\begin{bmatrix}
        1 & 0\\
        0      & v^2+w^2
       \end{bmatrix} =\begin{bmatrix}
        \rho^2 & 0\\
        0      & \rho^2\dot{f}^2 + \dot{\rho}^2
       \end{bmatrix}.
    \end{align*}

We note that, since $g_{\beta \beta} = \dot{\rho}^2 + \rho^2\dot{f}^2$, for $F$ to be an immersion, we cannot have $\dot{\rho}$ and $\dot{f}$ vanishing at the same time. This is equivalent to $\tau$ being nowhere zero, which is, in turn, equivalent to $\gamma(\beta) = \rho(\beta)e^{if(\beta)}$ being a regular parametrization of $\gamma$ as $\dot{\gamma} = \tau \gamma$. Therefore, we will always assume that we are working with a regular parametrization of $\gamma$. 
    
The Christoffel symbols can be calculated by
 \begin{align*}
     \Gamma^k_{ij} = \frac{1}{2}g^{kl}\left(\partial_i g_{jl} + \partial_j g_{il} - \partial_l g_{ij}\right).
 \end{align*}
Therefore, 
 \begin{align*}
      \Gamma^\alpha_{\alpha \alpha} &= \frac{1}{2}g^{\alpha \alpha}\left(\partial_\alpha g_{\alpha \alpha} + \partial_\alpha g_{\alpha \alpha} - \partial_\alpha g_{\alpha \alpha}\right) =  \frac{1}{2}g^{\alpha \alpha}\partial_\alpha g_{\alpha \alpha}= \frac{1}{2\rho^2} \partial_\alpha \big( \rho^2 \big)= 0\\
      \Gamma^\alpha_{\alpha \beta} &=  \Gamma^\alpha_{\beta \alpha} =  \frac{1}{2}g^{\alpha \alpha}\left(\partial_\beta g_{\alpha \alpha} + \partial_\alpha g_{\beta \alpha} - \partial_\alpha g_{\beta \alpha}\right) =  \frac{1}{2}g^{\alpha \alpha} \partial_\beta g_{\alpha \alpha}=  \frac{1}{2\rho^2} \partial_\beta \big( \rho^2 \big)  = \frac{\dot{\rho}}{\rho} = v\\
      \Gamma^\alpha_{\beta \beta} &= \frac{1}{2}g^{\alpha \alpha}\left(\partial_\beta g_{\beta \alpha} + \partial_\beta g_{\beta \alpha} - \partial_\alpha g_{\beta \beta}\right) =  -\frac{1}{2}g^{\alpha \alpha}\partial_\alpha g_{\beta \beta}= -\frac{1}{2\rho^2} \partial_\alpha \big( \rho^2\big(v^2 +w^2\big) \big)= 0\\
       \Gamma^\beta_{\beta \beta} &= \frac{1}{2}g^{\beta \beta}\left(\partial_\beta g_{\beta \beta} + \partial_\beta g_{\beta \beta} - \partial_\beta g_{\beta \beta}\right) =  \frac{1}{2}g^{\beta \beta}\partial_\beta g_{\beta \beta}=\frac{1}{2\rho^2\left(v^2 +w^2\right)} \partial_\beta \big(\rho^2\big(v^2 +w^2\big)\big) = v +  \frac{v\dot{v} + w\dot{w}}{v^2 + w^2}\\
  \Gamma^\beta_{\alpha \beta} &=  \Gamma^\beta_{\beta \alpha} =  \frac{1}{2}g^{\beta \beta}\left(\partial_\beta g_{\alpha \beta} + \partial_\alpha g_{\beta \beta} - \partial_\beta g_{\beta \alpha}\right) =  \frac{1}{2}g^{\beta \beta} \partial_\alpha g_{\beta \beta}=\frac{1}{2 \rho^2\left(v^2 +w^2\right)} \partial_\alpha \big( \rho^2\big(v^2 +w^2\big)\big) = 0\\
  \Gamma^\beta_{\alpha \alpha} &= \frac{1}{2}g^{\beta \beta}\left(\partial_\alpha g_{\alpha \beta} + \partial_\alpha g_{\alpha \beta} - \partial_\beta g_{\alpha \alpha}\right) =  -\frac{1}{2}g^{\beta \beta}\partial_\beta g_{\alpha \alpha}= -\frac{1}{2 \rho^2\left(v^2 +w^2\right)} \partial_\alpha \big(\rho^2 \big)=- \frac{v}{v^2 + w^2}.
 \end{align*}

Let $B$ denote the second fundamental form of $L_\gamma$ in $\C^2$. Since $\C^2$ is flat, $B$ is given by the covariant Hessian of $F$. Let $\{e^1, e^2\}$ be the dual basis corresponding to $\{e_1, e_2\}$ and let $(\Tilde{x}^1,\Tilde{x}^2,\Tilde{x}^3,\Tilde{x}^4) =(x^1,y^1,x^2,y^2)$ denote the usual Euclidean coordinate functions on $\R^4 \cong \C^2$. Then the components of 
 \begin{align*}
     B = B^A_{ij} e^i \otimes e^j \otimes \frac{\partial}{\partial \Tilde{x}^A}
 \end{align*}
 are given by 
 \begin{align*}
     B^A_{ij} = {\rm Hess}_g(F^A)_{ij} = \partial^2_{ij} F^A -  \Gamma^k_{ij}\partial_k F^A
 \end{align*}
 where we abused notation by using $\partial_1$ and $\partial_2$ to denote $\partial_\alpha$ and $\partial_\beta$ respectively. So we have
 \begin{align*}
      B_{\alpha \alpha} &= \partial^2_{\alpha \alpha} F - \Gamma^k_{\alpha \alpha}\partial_k F \\
      &= \partial^2_{\alpha \alpha} F - \Gamma^\beta_{\alpha \alpha}\partial_\beta F \\
      &=-F + \frac{v}{v^2 + w^2} \tau F\\
      &= \frac{iw( v +i w)}{v^2 + w^2}F\\
      &= \frac{w}{v^2 + w^2}Je_2\\
       B_{\beta \alpha} &=  B_{\alpha \beta } = \partial^2_{\alpha \beta } F - \Gamma^k_{ \alpha \beta }\partial_k F\\
      &=\partial^2_{ \alpha \beta } F - \Gamma^\alpha_{ \alpha \beta }\partial_\alpha F\\
      &=\partial_{\alpha} \left( \tau F\right) -v\partial_{\alpha}F \\
      &= \left( \tau -v \right)\partial_{\alpha}F \\
      &= iw\partial_{\alpha}F \\
      &= w Je_1\\
%\text{and}\\
      B_{\beta \beta} & = \partial^2_{\beta \beta} F - \Gamma^k_{\beta \beta}\partial_k F \\
      &=\partial^2_{\beta \beta} F - \Gamma^\beta_{\beta \beta}\partial_\beta F\\ 
      &=\partial_{\beta} \left( \tau F\right) -\left( v +  \frac{w \dot{w} + v\dot{v}}{v^2 +w^2} \right) \tau F \\
      &= \dot{\tau}F + \tau^2F -\left( v +  \frac{w \dot{w} + v\dot{v}}{v^2 +w^2} \right)  \tau F\\
      &= \left(\frac{\dot{\tau}}{\tau} + \tau -v -  \frac{w \dot{w} + v\dot{v}}{v^2 +w^2} \right) \tau F \\
      &= \left( \frac{(\dot{v} + i\dot{w})(v-iw)}{v^2 +w^2}+ iw - \frac{w \dot{w} + v\dot{v}}{v^2 +w^2} \right) \tau F\\
       &=i \left(w + \frac{ \dot{w}v - w\dot{v}}{v^2 +w^2} \right) \tau F \\
       &=  \left(w + \frac{ \dot{w}v - w\dot{v}}{v^2 +w^2} \right) Je_2.
\end{align*}
Now, we can calculate the mean curvature vector $H$ of $L_\gamma$.
\begin{lemma}\label{lem:H_C_rho,f}
     The mean curvature of $L_\gamma$ is given by
     \[
     H = C_{\rho,f}Je_2
     \]
     where 
     \[
        C_{\rho,f} =\frac{\dot{w}v-w\dot{v} +2w\left(v^2 +w^2\right)}{\rho^2\left(v^2+w^2\right)^2} =\frac{\rho\dot{\rho}\ddot{f}+2\rho^2\dot{f}^3+3\dot{\rho}^2\dot{f}- \rho \ddot{\rho}\dot{f} }{(\rho^2\dot{f}^2 + \dot{\rho}^2)^2}.
     \]
 \end{lemma}

 \begin{proof}
     The formula follows from straightforward calculations. Taking the trace of the second fundamental form $B$ gives us
     \begin{align*}
         H &= g^{\alpha \alpha} B_{\alpha\alpha} + g^{\beta \beta} B_{\beta\beta}\\
         &=\frac{1}{\rho^2}\frac{w}{v^2 + w^2}Je_2 + \frac{1}{\rho^2\left(v^2 +w^2\right)}\left(w + \frac{ \dot{w}v - w\dot{v}}{v^2 +w^2} \right) Je_2\\
         &= \frac{\dot{w}v-w\dot{v} +2w\left(v^2 +w^2\right)}{\rho^2\left(v^2+w^2\right)^2} Je_2
     \end{align*}
 \end{proof}
 
\begin{remark}
    Let us define the orthonormal basis $$\epsilon_1 = \frac{e_2}{\abs{e_2}}=\frac{e_2}{\rho \sqrt{v^2+ w^2}}, \,\, \epsilon_2 = \frac{e_1}{\abs{e_1}} =\frac{e_1}{\rho}.$$ In terms of this orthonormal basis,
\begin{align*}
    B(\epsilon_1,\epsilon_1) &= \frac{1}{\rho^2 (v^2+ w^2)}B_{\beta \beta} \\
    &= \frac{1}{\rho^2 (v^2+ w^2)} \left( \left(w + \frac{ \dot{w}v - w\dot{v}}{v^2 +w^2} \right) Je_2 \right)\\ &=\left(\frac{w}{\rho\sqrt{w^2 +v^2}} + \frac{ \dot{w}v - w\dot{v}}{\rho\left(v^2 +w^2\right)^{\frac{3}{2}}} \right) J\epsilon_1 \\
    B(\epsilon_2,\epsilon_2) &=  \frac{1}{\rho^2}B_{\alpha \alpha}= \frac{1}{\rho^2} \left(\frac{w}{w^2 +v^2} Je_2 \right) = \frac{w}{\rho \sqrt{w^2 +v^2}} J\epsilon_1 \\
     B(\epsilon_1,\epsilon_2) &=\frac{1}{\rho^2\sqrt{w^2 +v^2}}B_{\beta \alpha} = \frac{1}{\rho^2\sqrt{w^2 +v^2}}(w Je_1) = \frac{w}{\rho\sqrt{w^2 +v^2}}J\epsilon_2.
\end{align*}
Therefore, the immersion is H-umbilical according to \citep{BY-Chen_H-umbilical}. 
\end{remark}

Next, we calculate the divergence of the vector field $JH$ and prove an identity that is useful for characterizing Hamiltonian stationary Lagrangian twisted tori.
\begin{lemma}\label{lemma:div_JHs}
The divergence of $JH$ is given by
\begin{align}\label{eqn:divJH_C_rho,f}
     {\rm div}_g(JH) &=  -\left[ \frac{d}{d \beta} C_{\rho,f} + \left(\frac{d}{d \beta} \ln{ \sqrt{\det g}} \right)C_{\rho,f}\right].
\end{align}
   Furthermore, the mean curvature vector $H$ of $L_\gamma$ satisfies
 \[
    \nabla \left( \rho^2 \abs{JH}^2\right) = 2\rho^2 {\rm div}_g(JH) JH \numberthis \label{eqn:grad_rho2_JH2}.
 \]
 
\end{lemma}
\begin{proof}
    By Lemma \ref{lem:H_C_rho,f}, $JH = -C_{\rho,f}e_2$. Therefore,
 \begin{align*}
     \Div_g(JH) &= g^{\alpha \alpha} \left\langle \partial_\alpha JH, e_1 \right\rangle_{Euc} + g^{\beta \beta} \left\langle \partial_\beta JH, e_2 \right\rangle_{Euc}\\
     &=-g^{\alpha \alpha} C_{\rho,f}\Gamma^{\alpha}_{\alpha \beta}g_{\alpha \alpha} - g^{\beta \beta}\frac{d}{d \beta} C_{\rho,f}g_{\beta \beta} - g^{\beta \beta} C_{\rho,f}\Gamma^{\beta}_{\beta \beta}g_{\beta \beta}\\
     &= -\left[ \frac{d}{d \beta} C_{\rho,f} + \left(\Gamma^{\alpha}_{\alpha \beta} + \Gamma^{\beta}_{\beta \beta}\right) C_{\rho,f}\right]\\
     &=  -\left[ \frac{d}{d \beta} C_{\rho,f} + \left(\frac{d}{d \beta} \ln{ \sqrt{\det g}} \right)C_{\rho,f}\right]
 \end{align*}
where, in the last equality, we used the identity $\sum_i \Gamma^{i}_{ij} = \partial_{j}\ln{\sqrt{\det{g}}}$. This proves (\ref{eqn:divJH_C_rho,f}). 

Also, the right-hand side of (\ref{eqn:grad_rho2_JH2}) can be written as
\begin{align*}
    \rho^2 \Div_g(JH) JH &= \rho^2\left[ \frac{d}{d \beta} C_{\rho,f} + \left(\frac{d}{d \beta} \ln{ \sqrt{\det{g}}} \right)C_{\rho,f}\right]C_{\rho,f} e_2\\
                        &=\left[\rho^2\frac{d}{d \beta} C_{\rho,f}C_{\rho,f} + \rho^2\left(\frac{d}{d \beta} \ln{ \sqrt{\det{g}}} \right)C_{\rho,f}^2\right] e_2\\
                        &=\left[\rho^2\frac{d}{d \beta} C_{\rho,f}C_{\rho,f} + \frac{1}{2}g_{\alpha \alpha}\left(\frac{d}{d \beta} \ln{( g_{\alpha \alpha} g_{\beta \beta})} \right)C_{\rho,f}^2\right] e_2\\
                        &=\left[\rho^2\frac{d}{d \beta} C_{\rho,f}C_{\rho,f} + \frac{1}{2} g^{\beta \beta}\left(\frac{d}{d \beta} ( g_{\alpha \alpha} g_{\beta \beta}) \right)C_{\rho,f}^2\right] e_2.
\end{align*}
On the other hand, $ \rho^2 \abs{JH}^2 = \rho^2g_{\beta \beta}C_{\rho,f}^2$ so the left-hand side of (\ref{eqn:grad_rho2_JH2}) can be written as
\begin{align*}
    \nabla \left( \rho^2 \abs{JH}^2\right) &= g^{\beta \beta}\frac{d}{d \beta}\left(\rho^2 \abs{JH}^2\right)e_2\\
                                            &= g^{\beta \beta}\frac{d}{d \beta}\left( \rho^2g_{\beta \beta}C_{\rho,f}^2\right)e_2\\
                                            &=\left[g^{\beta \beta}\frac{d}{d \beta}\left( \rho^2g_{\beta \beta}\right)C_{\rho,f}^2+ 2\rho^2C_{\rho,f} \frac{d}{d \beta} C_{\rho,f} \right]e_2\\
                                            &=\left[g^{\beta \beta}\frac{d}{d \beta}\left( g_{\alpha \alpha} g_{\beta \beta}\right)C_{\rho,f}^2+ 2\rho^2C_{\rho,f} \frac{d}{d \beta} C_{\rho,f} \right]e_2
\end{align*}
which proves (\ref{eqn:grad_rho2_JH2}).
\end{proof}
We conclude this section by taking a look at how the geometry of $\gamma$ in $\C$ relates to the geometry of $L_\gamma$ in $\C^2$.
\begin{lemma}\label{lem:gamma_gometry}
Let $\gamma$ be an embedded simple closed curve in $\C^2$ parametrized by $\gamma(\beta) = \rho(\beta)e^{if(\beta)}$ where $\rho:\R \to \R$ is a positive $2\pi$-periodic function and $f: \R \to \R$ is a function that satisfies $f(x+2\pi) = f(x) +2k\pi$ for some $k \in \Z$. Then
\begin{enumerate}
    \item the winding number of $\gamma$ around the origin is $k$ and therefore $f$ is a $2\pi$-periodic function when $L_\gamma$ is a Chekanov torus.
    %%\item $L_\gamma$ is an exotic torus if and only if $f$ is a $2\pi$-periodic function;
    \item $L_\gamma$ is the product of two circles if and only if $\rho$ is constant.
\end{enumerate}
\end{lemma}
\begin{proof}
\begin{enumerate}
    \item    The winding number of $\gamma$ around the origin is given by
\begin{align*}
        \textbf{ind}(\gamma,0) = \frac{1}{2\pi i} \int_\gamma \frac{dz}{z} = \frac{1}{2\pi i} \int_0^{2 \pi} \frac{\dot{\gamma}}{\gamma} dt.
 \end{align*}
Since $\dot{\gamma}  =\big(\dot{\rho} + i\rho \dot{f}\big)e^{if} = \tau \gamma = (v+ iw)\gamma$ and $\rho$ is a $2\pi$-periodic function,
\begin{align*}
        \textbf{ind}(\gamma,0) &= \frac{1}{2\pi i} \left(\int_0^{2 \pi} v dt + i\int_0^{2 \pi} w dt\right)\\
        &= \frac{1}{2\pi i} \left(\int_0^{2 \pi} \frac{d }{dt}\ln{\rho} dt + i\int_0^{2 \pi} \dot{f} dt\right) \\
        &= \frac{1}{2 \pi}\left( f(2\pi) - f(0)\right)\\
        &=k.
\end{align*}
When $L_\gamma$ is one of Chekanov's exotic tori then $\gamma$ is contained in an open half-plane and hence its winding number around the origin is $0$. Therefore, $k=0$ and $f$ is a $2\pi$-periodic function, which finishes the proof of the first statement.
\item Suppose that $\rho$ is constant. After performing the change of coordinates $x=\beta + \alpha$ and $y=\beta - \alpha$ for $\alpha, \beta \in [0,2\pi)$, the immersion $F$ takes the form
\[
    F(x,y) = \frac{\rho}{\sqrt{2}}\left(e^{ix},e^{iy}\right)
\]
where $(x,y) \in D =\big\{(x,y) \in \R^2 \big| -x \leq y < 4\pi-x \text{ and } x-4\pi < y \leq x \big\}$. In these coordinates, it is clear that the image of $F$ is the product torus $S^1\left(\frac{\rho}{\sqrt{2}}\right) \times S^1\left(\frac{\rho}{\sqrt{2}}\right) \subset \C^2$ where, for a positive constant $r$, $S^1\left( r\right):=\left\{z\in \C \big| \,\abs{z} = r \right\}$. On the other hand, if $L_\gamma$ is the product $S^1\left(r_1\right) \times S^1\left( r_2\right)$, then $\frac{\rho}{\sqrt{2}} = \abs{\frac{\rho}{\sqrt{2}}e^{i\left(f+\alpha\right)}} = r_1$ and  $\frac{\rho}{\sqrt{2}} = \abs{\frac{\rho}{\sqrt{2}}e^{i\left(f-\alpha\right)}} = r_2$ so we must have that $\rho = \sqrt{2}r_1 = \sqrt{2}r_2$ is a constant.
\end{enumerate}
\end{proof}

 \section{Minimality of the twisted torus $L_\gamma$}\label{section:minimality}
We start this section by reformulating the Hamiltonian stationary equation (\ref{eqn:E-L_general}) in terms of the functions $v$ and $w$.
\begin{theorem}\label{thm:E-L_gamma}
The twisted torus $L_\gamma$ is Hamiltonian stationary, if and only if
 \[
    \sqrt{\det{g}} C_{\rho,f} =\frac{\dot{w}v-w\dot{v} +2w\left(v^2 +w^2\right)}{\left(v^2+w^2\right)^{\frac{3}{2}}} = c \numberthis \label{eqn:E-L_gamma}
 \] 
 for some non-zero constant $c$. Moreover, a Hamiltonian stationary Lagrangian twisted torus satisfies
 \[
    \frac{d}{d \beta}\left(\rho^2 \left(\phi -\frac{c}{2} \right)\right)=0 \numberthis \label{eqn:E-L_phi}
 \]
 where $\phi = \frac{w}{\sqrt{v^2 +w^2}}$.
 \end{theorem}

 \begin{proof}
  By (\ref{eqn:divJH_C_rho,f}) of Lemma \ref{lemma:div_JHs}, $L_\gamma$ is Hamiltonian stationary if and only if 
\begin{align}
    \frac{d}{d \beta} C_{\rho,f} + \left(\frac{d}{d \beta} \ln{ \sqrt{\det{g}}} \right)C_{\rho,f} = 0,
\end{align}
or equivalently,
\begin{align*}
   \sqrt{\det{g}} \frac{d}{d \beta} C_{\rho,f} + \left(\frac{d}{d \beta}  \sqrt{\det{g}} \right)C_{\rho,f} = 0 \\
   \frac{d}{d \beta} \left( \sqrt{\det{g}} C_{\rho,f}\right) = 0.
\end{align*}
Therefore, $\sqrt{\det{g}} C_{\rho,f} = c$ for some constant $c$. Expressing $\sqrt{\det{g}}$ and $C_{\rho,f}$ in terms of $v$ and $w$ gives us (\ref{eqn:E-L_gamma}). Note that $c=0$ if and only if $H = 0$ by Lemma \ref{lem:H_C_rho,f} but there are no minimal tori in $\C^2$ so $c$ cannot be zero.

For  $\phi = \frac{w}{\sqrt{v^2 +w^2}}$, we have 
\begin{align*}
    \dot{\phi}&=v\frac{\dot{w}v -w\dot{v}}{\left( w^2 + v^2\right)^{\frac{3}{2}}}
\end{align*}
so multiplying (\ref{eqn:E-L_gamma}) by $v$ gives us
\begin{align*}
    \dot{\phi} + 2v\phi = cv, 
\end{align*} or equivalently
\begin{align*}
    \dot{\phi} + 2v\left(\phi -\frac{c}{2} \right) = 0.
\end{align*}
Finally, after multiplying both sides by $\rho^2$, we get that
\begin{align*}
    0 &= \rho^2\dot{\phi} + 2\rho \dot{\rho}\left(\phi -\frac{c}{2}  \right) \\
    &=\rho^2\frac{d}{d \beta}\left(\phi -\frac{c}{2}  \right)  + 2\rho \dot{\rho}\left(\phi -\frac{c}{2} \right)\\
    &=\frac{d}{d \beta}\left(\rho^2 \bigg(\phi -\frac{c}{2}  \bigg)\right)
\end{align*}
which proves (\ref{eqn:E-L_phi}).
 \end{proof}

Next, we state an immediate corollary of Lemma \ref{lemma:div_JHs}, which characterizes Hamiltonian stationary Lagrangian twisted tori.
\begin{cor}\label{cor:constant_rho_JH}
    The twisted torus $L_\gamma$ is Hamiltonian stationary if and only if $\rho \abs{H}$ is constant.
\end{cor}
\begin{proof}
    If $L_\gamma$ is Hamiltonian stationary, then $\Div_g(JH)=0$ so, by (\ref{eqn:grad_rho2_JH2}), $\nabla \left( \rho^2 \abs{JH}^2\right)=0$.
    On the other hand, if $ \nabla \left( \rho^2 \abs{JH}^2\right) = 0$, then by (\ref{eqn:grad_rho2_JH2}), we must have $\Div_g(JH)(x)=0$ or $JH(x) = 0$ at each $x\in L_\gamma$. Fix $x_0 \in L_\gamma$. If  $\Div_g(JH)(x_0)=0$, then we are done, so let us assume that $JH(x_0) = 0$. If $x_0$ has a neighbourhood in which $JH =0$, then clearly $\Div_g(JH)(x_0)=0$. Finally, if $x_0$ has no such neighbourhood, then there exist $x_n \to x_0$ such that $\Div_g(JH)(x_n)=0$, so we must again have $\Div_g(JH)(x_0)=0$ by the continuity of $\Div_g(JH)$. Therefore, we can conclude that we must have $\Div_g(JH)(x_0)=0$ and, since $x_0$ was chosen arbitrarily, that $\Div_g(JH)=0$ everywhere.
\end{proof}

 We also state two useful corollaries of Theorem \ref{thm:E-L_gamma}. Corollary \ref{cor:product_torus} characterizes the case when a generic twisted torus is the product of two circles, while Corollary \ref{cor:star_shaped} gives us a sufficient condition for a Hamiltonian stationary twisted torus to be a product.
\begin{cor}\label{cor:product_torus}
    The torus $L_\gamma$ is a product torus if and only if $\rho \abs{H} = 2$.
\end{cor}
\begin{proof}
     By Lemma \ref{lem:gamma_gometry}, $L_\gamma$ is a product torus if and only if $\rho$ is constant, so it is sufficient to show that $\rho\abs{H}=2$ if and only if $\rho$ is constant.
     
    First, suppose that $\rho \abs{H}=2$. Corollary \ref{cor:constant_rho_JH} tells us that $L_\gamma$ is Hamiltonian stationary, so by (\ref{eqn:E-L_gamma}) of Theorem \ref{thm:E-L_gamma},
    \[
    \sqrt{\det{g}} C_{\rho,f} = c  \label{eqn:detg_C}
    \]
    for some non-zero constant $c$. Since $\sqrt{\det{g}} = \sqrt{g_{\alpha \alpha}g_{\beta \beta}}$ and, by Lemma \ref{lem:H_C_rho,f}, 
    $$
    \rho^2\abs{JH}^2 = \rho^2 g_{\beta \beta}C_{\rho,f}^2 = g_{\alpha \alpha}g_{\beta \beta}C_{\rho,f}^2,
    $$
     we have
    \[
    \rho^2 \abs{JH}^2 = c^2 \numberthis \label{eqn:rho2_JH2}
    \]
    and thus $\abs{c} = 2$. Also, by (\ref{eqn:E-L_phi}) of Theorem \ref{thm:E-L_gamma},
    \begin{align*}
        \frac{d}{d \beta}\left(\rho^2 \left(\phi -\frac{c}{2}  \right)\right) = 0
    \end{align*}
or equivalently,
\begin{align}\label{eqn:E-L_phi_2}
    \rho^2 \left(\phi -\frac{c}{2} \right) = b
\end{align}
for some constant $b$. Since $\rho$ attains both its global minimum and global maximum along $\gamma$, it has at least two distinct critical points in $[0,2 \pi)$. Let $\beta_{max}$ and $\beta_{min}$ be such that $\rho(\beta_{max}) = \rho_{max} := \max_{\gamma}\rho$ and $\rho(\beta_{min}) = \rho_{min} := \min_{\gamma}\rho$. Then
\begin{align*}
  \rho_{max}^2\left(\phi(\beta_{max}) -\frac{c}{2} \right)  = b = \rho_{min}^2\left(\phi(\beta_{min}) -\frac{c}{2} \right).
\end{align*}
Note that, since $\phi=\frac{w}{\sqrt{v^2 +w^2}}$ by definition and $v(\beta_{max}) = v(\beta_{min}) = 0$, we have $\abs{\phi(\beta_{max})} = \abs{\phi(\beta_{min})} = 1 = \frac{\abs{c}}{2}$.

If $b \neq 0$, then we must have $\phi(\beta_{max}) = \phi(\beta_{min}) = -\frac{c}{2}$ and thus
\begin{align*}
 -c \rho_{max}^2  = -c\rho_{min}^2.
\end{align*}
Therefore, $ \rho_{max} =\rho_{min}$ and we can conclude that $\rho$ is constant. However, when $\rho$ is constant, then both $v$ and $\dot{v}$ are identically zero, so
\begin{align*}
    \phi =  \frac{w}{\abs{w}}=\frac{w^3}{\abs{w}^3} = \frac{c}{2}
\end{align*}
by (\ref{eqn:E-L_gamma}) and thus $b=0$ by (\ref{eqn:E-L_phi_2}). This is a contradiction so we can conclude that $b=0$. Then, ${|\phi| = \frac{\abs{c}}{2}} = 1$ identically by (\ref{eqn:E-L_phi_2}) which in turn implies that $v$ is identically zero. Therefore, $\rho$ must be constant.

\begin{comment}
  On the other hand, when $\rho$ is constant, then both $v$ and $\dot{v}$ are zero. So the left-hand side of equation (\ref{eqn:E-L_gamma}) simplifies to
    \[
        \frac{2w^3}{(w^2 )^{\frac{3}{2}}} = 2\frac{\dot{f}^3}{\abs{\dot{f}}^3} = -2 \text{sign}\big( \dot{f} \big)^3.
    \]
    We assume that $\gamma$ has a regular parametrization so $\dot{\rho}$ and $\dot{f}$ cannot vanish at the same time. Since $\rho$ is constant, $\dot{f}$ can never vanish and hence it can never change sign. So $\text{sign}\big( \dot{f} \big)$ is constant and we can conclude that $L_\gamma$ is Hamiltonian stationary with $\abs{c} = 2$. Therefore, by equation (\ref{eqn:rho2_JH2}), $\rho \abs{H}=2$.
\end{comment}

On the other hand, when $\rho$ is constant, then $L_\gamma$ is the product torus $S^1\left(\frac{\rho}{\sqrt{2}}\right) \times S^1\left(\frac{\rho}{\sqrt{2}}\right)$ by Lemma \ref{lem:gamma_gometry}. Therefore, $\rho \abs{H} = 2$ by direct calculations.
\end{proof}

\begin{cor}\label{cor:star_shaped}
      Let $\gamma: S^1 \to \C$ be a simple closed curve that encloses a star-shaped region centred at the origin. Then $L_\gamma$ is Hamiltonian stationary if and only if it is the product torus $S^1(r)\times S^1(r)$ where $2\pi r^2$ equals the area enclosed by $\gamma$.  
 \end{cor}
 \begin{proof}
    Let $\gamma: S^1 \to \C$ be a simple closed curve that encloses a star-shaped region centred at the origin. Then $\gamma$ can be parametrized as $\gamma(\beta) = \rho(\beta)e^{i\beta}$ for some positive $2 \pi$-periodic function $\rho$, i.e. we may take $f(\beta) = \beta$, and $\phi = \frac{1}{\sqrt{1 + v^2}}$.
    
    First, let us assume that $L_\gamma$ is Hamiltonian stationary. Then,  by (\ref{eqn:E-L_phi}) of Theorem \ref{thm:E-L_gamma},
    \begin{align*}
        \frac{d}{d \beta}\left(\rho^2 \left(\phi -\frac{c}{2} \right)\right) = 0
    \end{align*}
or equivalently,
\[
    \rho^2 \left(\phi  -\frac{c}{2} \right) = b
\]
for some constant $b$. Therefore,
\begin{align*}
  \rho_{max}^2\left(\phi(\beta_{max}) -\frac{c}{2} \right)  = b = \rho_{min}^2\left(\phi(\beta_{min}) - \frac{c}{2} \right)
\end{align*}
where $\beta_{min},\beta_{max},\rho_{min}$ and $\rho_{max}$ are defined as in the proof of Corollary \ref{cor:product_torus}. We know that $v(\beta_{max}) = v(\beta_{min})= 0$ so $\phi(\beta_{max}) = \phi(\beta_{min}) = 1$ and
\[
\rho_{max}^2\left(1 -\frac{c}{2} \right)  =\rho_{min}^2\left(1 - \frac{c}{2} \right).
\]
 If $c\neq 2$, then we must have $\rho_{min} = \rho_{max}$, so $\rho$ is constant and $c=2$ by Theorem \ref{thm:E-L_gamma}. This is a contradiction so we can conclude that $c$ must be $2$. Then, since $\abs{c} = \rho\abs{JH}$ by (\ref{eqn:rho2_JH2}), $L_\gamma$ is a product torus by Corollary \ref{cor:product_torus}.

 The converse statement follows from Corollary \ref{cor:product_torus} and Corollary \ref{cor:constant_rho_JH}.
 \end{proof}

We finish this section by stating two lemmas about Hamiltonian stationary twisted tori that will help us to prove Theorem \ref{theorem:main}.
\begin{lemma}\label{lemma:curvature}
    Suppose that $\gamma$ is a simple closed curve in $\C$ that is oriented counterclockwise. If $L_\gamma$ is Hamiltonian stationary, then the signed curvature of $\gamma$ is given by
    \begin{align*}
        \kappa = \frac{c-\phi}{\rho}
    \end{align*}
    where $c$ and $\phi$ are as defined in Theorem \ref{thm:E-L_gamma}.
\end{lemma}
\begin{proof}
    Suppose that $\gamma$ is parametrized counterclockwise. Then the inwards pointing unit normal vector to $\gamma$ is given by $N = \frac{J\dot{\gamma}}{\abs{\dot{\gamma}}}$ so the signed curvature of $\gamma$ can be calculated as
\begin{align*}
    \kappa &= \langle \nabla_{\partial s} T, N\rangle \\ 
    &=\text{Re}\left( \frac{d \beta}{ds}\frac{d}{d \beta}\left(\frac{\dot{\gamma}}{{|\dot{\gamma}|}}\right) \cdot \frac{\overline{(i \dot{\gamma})}}{{|\dot{\gamma}|}}\right)\\
    &=-\frac{1}{{|\dot{\gamma}|}^2}\text{Re}\left(i\frac{\ddot{\gamma} \cdot \overline{\dot{\gamma}}}{|{\dot{\gamma}}|} \right)
\end{align*}
where $s = \int_0^\beta \abs{\dot{\gamma}}d\beta$ is the arc-length parameter along $\gamma$. Recalling that $\dot{\gamma} = \tau \gamma$ for $\tau = v + iw$, we can write
\begin{align*}
    \kappa&=-\frac{1}{|{\dot{\gamma}}|^3}\text{Re}\left(i(\dot{\tau}\cdot \overline{\tau} |{\gamma}|^2+\tau |{\gamma}|^2\abs{\tau}^2)\right)
\end{align*}
which simplifies to
\begin{align*}
     \kappa &=\frac{1}{\rho}\left(\frac{\dot{w}v - w\dot{v}}{(v^2+w^2)^\frac{3}{2}} + \frac{w}{\sqrt{v^2 +w^2}}\right).
\end{align*}
Finally, by (\ref{eqn:E-L_gamma}) of Theorem \ref{thm:E-L_gamma},
\begin{align*}
    c = \frac{\dot{w}v - w\dot{v}}{(v^2+w^2)^\frac{3}{2}} + 2\frac{w}{\sqrt{v^2 +w^2}}
\end{align*}
so
\begin{align*}
    \kappa &=\frac{1}{\rho}\left(c - \frac{w}{\sqrt{v^2 +w^2}}\right) =\frac{c- \phi}{\rho}.
\end{align*}
\end{proof}

\begin{lemma}\label{lemma:crit_points}
    Let $L_\gamma$ be a Hamiltonian stationary twisted torus that arises from a simple closed curve $\gamma$ that does not enclose the origin. If $\rho = |{\gamma}|$ has exactly two critical values then it must also have exactly two critical points along $\gamma$.
\end{lemma}
\begin{proof}
Let $X$ be the set of critical points of $\rho$ along $\gamma$. First, we show that $X$ is a finite set. Assume, for a contradiction, that $X$ is infinite. Then $X$ has an accumulation point $x \in \gamma$. Since $X$ is closed, $x \in X$, i.e. $\dot{\rho}(x) = 0$, and since it can be approximated along $\gamma$ by a sequence $x_n \to x$ such that $\dot{\rho}(x_n)=0$, we must also have that $\ddot{\rho}(x) = 0$. Therefore, $\abs{c}=2$ by Lemma \ref{lem:H_C_rho,f} where $c$ is defined as in Theorem \ref{thm:E-L_gamma}. Recall that, by (\ref{eqn:rho2_JH2}), $ \rho \abs{JH} = \abs{c} =2$ so Corollary \ref{cor:product_torus} tells us that $L_\gamma$ is a product torus or equivalently that $\rho$ is constant. This contradicts the assumption that $\gamma$ does not enclose the origin so we can conclude that $X$ is a finite set.

Now, we prove that $X$ has exactly two elements. Since $\gamma$ is compact and $\rho$ is a non-constant smooth function over $\gamma$, we know that $\rho$ has at least two distinct critical points and at least two distinct critical values corresponding to the global maximum and global minimum of $\rho$ over $\gamma$.  Let $x_{min},x_{max} \in X$ be two critical points corresponding to a global minimum and a global maximum respectively. Assume for a contradiction, that there exists a third critical point $x\in X$. By our assumptions, $\rho$ has exactly two critical values so every critical point must be either a global maximum or a global minimum. Without the loss of generality, let us assume that $\Bar{x}_{max}:=x$ is a global maximum and that the three critical points can be ordered as $(x_{min},\Bar{x}_{max},x_{max})$ when tracing $\gamma$ clockwise. Notice that there must be another critical point $\Bar{x}_{min}$ corresponding to a global minimum somewhere between $\Bar{x}_{max}$ and $x_{max}$ otherwise there would be a constant arc of $\gamma$ joining these two points contradicting that $X$ is finite. Therefore, we have at least four critical points and they can be ordered as ($x_{min},\Bar{x}_{max},\Bar{x}_{min},x_{max}$) when tracing $\gamma$ clockwise. Let $\rho_{min} := \inf_\gamma \rho$ and $\rho_{max} := \sup_\gamma \rho$. Then $\gamma$ is contained in the annulus centred at the origin of inner radius $\rho_{min}$ and outer radius $\rho_{max}$. Also, $\gamma$ has at least two arcs connecting the inner and the outer circles of this annulus, namely, the arc joining $x_{min}$ and $\Bar{x}_{max}$ an the arc joining $\Bar{x}_{min}$ and $x_{max}$, see the left-hand side sketch of Figure \ref{fig:figure_1}.
\begin{figure}[h]
    \centering
    \includegraphics[width=0.98\linewidth]{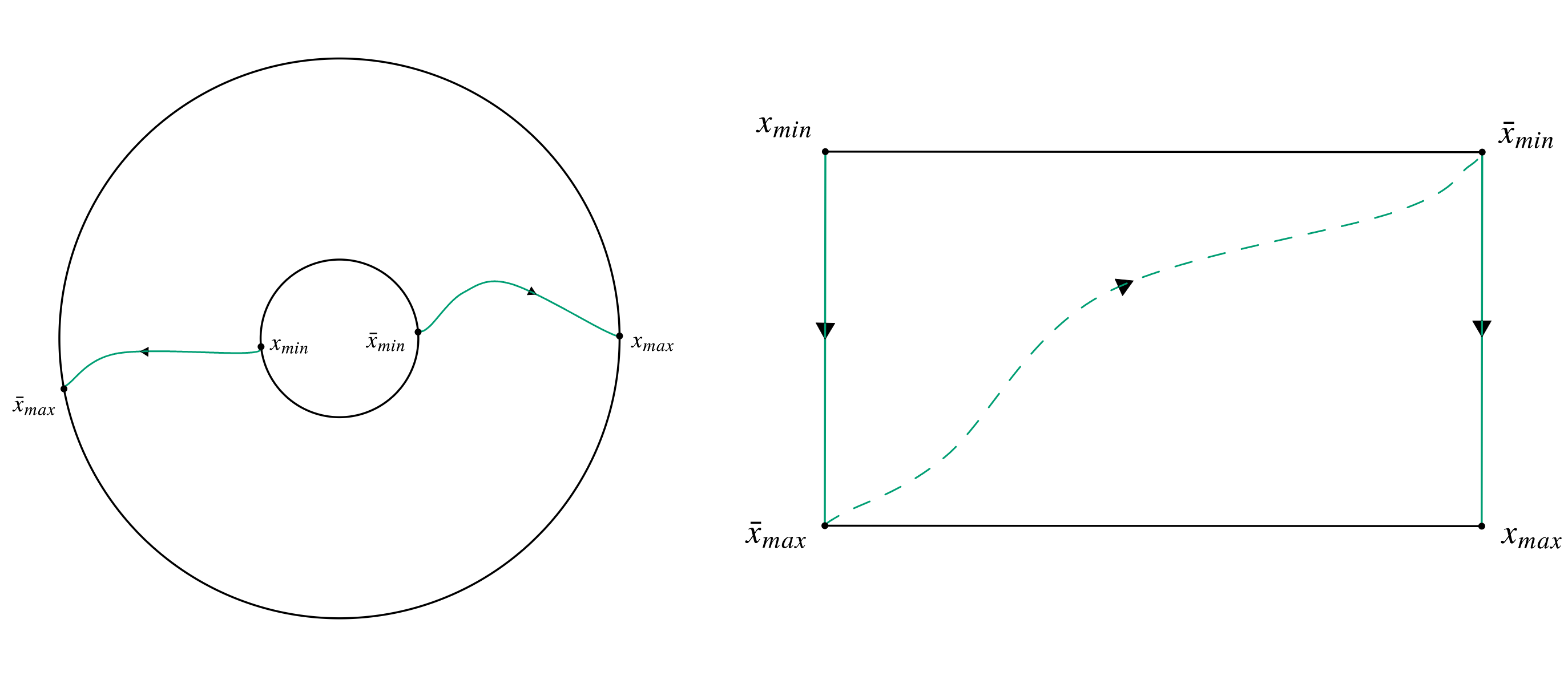}
    \caption{The sketch on the left shows the annulus of inner radius $\rho_{min}$ and outer radius $\rho_{max}$ together with the two arcs joining the two boundary circles. The sketch on the right shows the lower half of the annulus deformed into a rectangle.}
    \label{fig:figure_1}
\end{figure}
These two arcs divide the annulus into two halves, each of which can be visualized as rectangles, see the right-hand side sketch of Figure \ref{fig:figure_1}.  Since $\gamma$ is a closed curve, it must also have an arc connecting $\Bar{x}_{max}$ and $\Bar{x}_{min}$ and an arc connecting $x_{max}$ and $x_{min}$. Both of these arcs are connecting opposite vertices of one of the two halves of the annulus, see the right-hand side sketch of Figure \ref{fig:figure_1}. By our assumption, $\gamma$ has no self-intersection points, so these two arcs must be contained in different halves of the annulus as it is shown in Figure \ref{fig:figure_2}. This implies that $\gamma$ encloses the origin so we arrive at a contradiction and we can conclude that $X$ cannot have more than two elements.
\begin{figure}[h]
    \centering
    \includegraphics[width=0.5\linewidth]{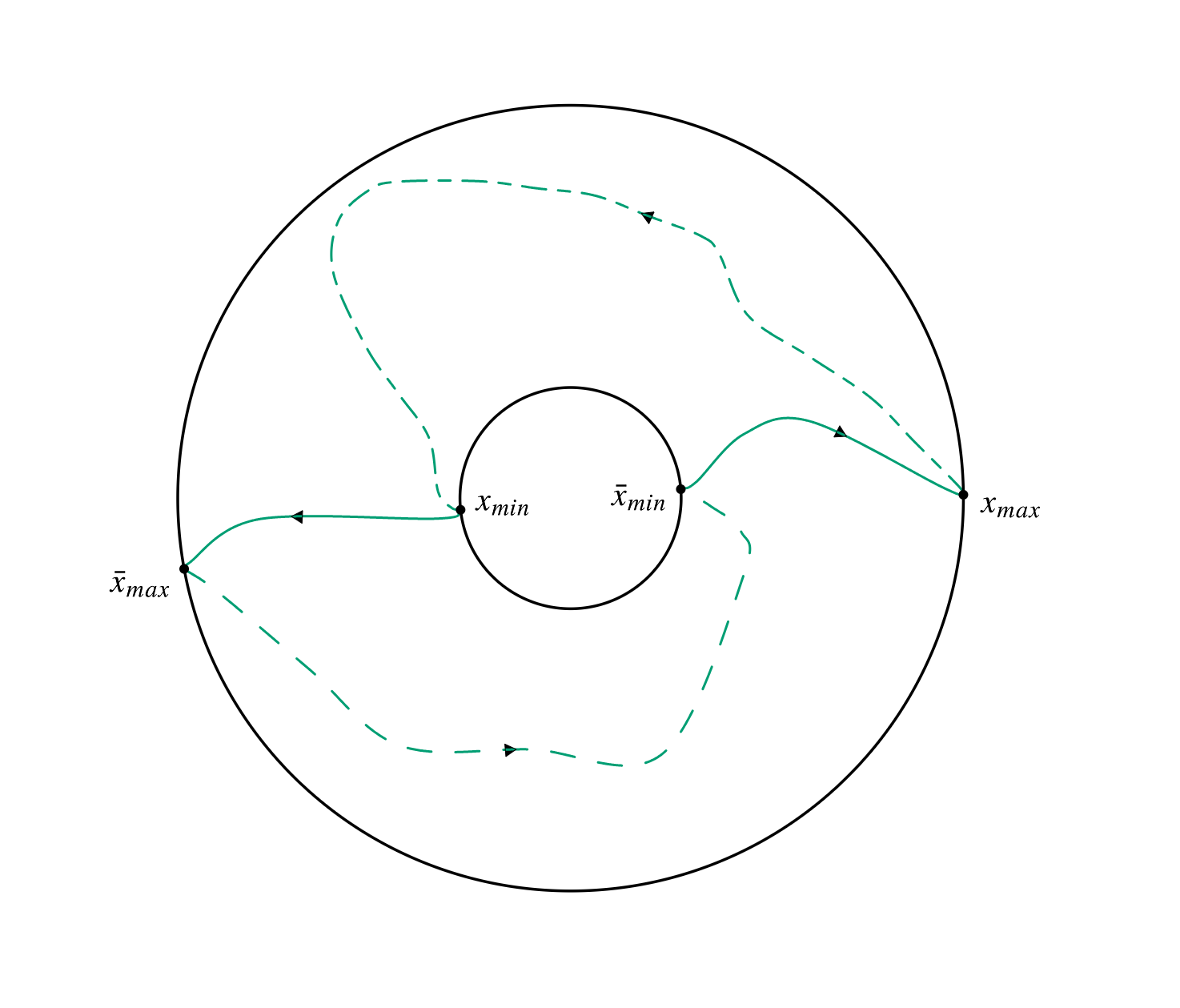}
    \caption{The only way $\gamma$ does not have a self-intersection point is if it encloses the origin. }
    \label{fig:figure_2}
\end{figure}
\end{proof}

\section{Proof of Theorem \ref{theorem:main} }\label{section:Proof}
Suppose that $L_\gamma$ is a twisted torus that corresponds to a simple closed curve $\gamma \subset  \C$ that is parametrized counterclockwise by $\gamma(\beta) = \rho(\beta)e^{if(\beta)}$ for $\beta \in [0,2 \pi)$. As discussed in Section \ref{section:introduction}., we may assume that $\rho$ is a positive, smooth, $2\pi$-periodic function and $f: \R \to \R$ is a smooth function that satisfies $f(x+2\pi) = f(x) +2k\pi$ for some $k \in \Z$.

First, assume that $\gamma$ is a circle centred at the origin. Then $\rho$ is constant so, by Lemma \ref{lem:gamma_gometry}, $L_\gamma$ is a product of two circles and hence it is Hamiltonian stationary by Corollary \ref{cor:product_torus} and Corollary \ref{cor:constant_rho_JH}.

Now, assume that  $L_\gamma$ is Hamiltonian stationary. We start by applying the change of variable $u(\beta) = \int_0^\beta \sqrt{v^2+w^2}$ where $w = \frac{df}{d \beta}$ and $v=\frac{d }{d \beta} \ln{\rho}= \frac{d \rho}{d \beta}\rho^{-1} $. Then
\begin{align*}
    \frac{du}{d\beta} = \sqrt{v^2+w^2}
\end{align*}
so
\[\frac{df}{du} =\frac{d\beta}{d u}\frac{df}{d\beta}  = \frac{w}{\sqrt{v^2 + w^2}} = \phi(\beta(u))\]

and
\[\frac{d \rho}{d u} \rho^{-1} = \frac{d\beta}{d u}\frac{d \rho}{d \beta} \rho^{-1} = \frac{v}{\sqrt{v^2 + w^2}}.\]
Therefore, changing to this new variable imposes the constraint
\[
\left(\frac{d \rho}{d u} \rho^{-1}  \right)^2 + \left(\frac{df}{du}\right)^2 = 1
\]
on $\rho$ and $f$.
Also, in terms of $u$, the Hamiltonian stationary equation (\ref{eqn:E-L_phi}) reads as 
\begin{align*}
    \frac{d}{du}\left(\rho^2\left( \frac{df}{du} -\frac{c}{2}\right)\right) = 0.
\end{align*}
Integrating the equation tells us that
\begin{align*}
    \frac{df}{du} = K\rho^{-2} +\frac{c}{2}
\end{align*}
for some constant $K$.
Integrating the equation a second time gives us
\begin{align*}
    f(u)  =\frac{c}{2}u + K I(u) + f(0)
\end{align*}
where  
$$
I(u) = \int_0^u \rho^{-2}(x)dx.
$$ 
Let $u^* = u(2\pi)$ and $I^* = I(u^*)$. We know that
\begin{align*}
    2k \pi = f(u^*) - f(0) =\frac{cu^*}{2} +KI^*,
\end{align*}
so 
$$
K =\frac{2k\pi}{I^*} -\frac{cu^*}{2I^*}
$$ 
and
\begin{align*}
    f(u) = \frac{cu^*}{2}\left(\frac{u}{u^*}-\frac{I(u)}{I^*}\right) +\frac{2k \pi }{I^*}I(u) +f(0)
\end{align*}

By Lemma \ref{lem:gamma_gometry}, $k$ is the winding number of the curve $\gamma$ around the origin. Since $\gamma$ is a simple closed curve that is parametrized counterclockwise, the only possible values for its winding number are $k=0$ and $k=1$. By Lemma \ref{lemma:curvature}, the signed curvature $\kappa$ of $\gamma$ is given by
\begin{align*}
    \kappa = \frac{c-\phi}{\rho}.
\end{align*}
Therefore, we have
\begin{align*}
    2  \pi &= \int_\gamma \kappa ds\\
    &= \int_0^{2\pi} \kappa \rho\sqrt{v^2 + w^2}d\beta\\
    &= \int_0^{2\pi} \frac{c- \phi}{\rho} \rho\sqrt{v^2 + w^2}d\beta\\ 
    &=c\int_0^{2\pi} \sqrt{v^2 + w^2}d\beta - \int_0^{2\pi} wd\beta\\
    &=cu^* -f(2\pi) + f(0)\\
    &=cu^* - 2k\pi
\end{align*}
which tells us that 
\begin{align}\label{eqn:c_u_star}
    u^* =\frac{2(k+1) \pi}{c}
\end{align}
and
\begin{align*}
    f(u) = \frac{(k+1) \pi}{u^*}u + \frac{(k-1) \pi}{I^*}I(u) +f(0).
\end{align*}

If $k = 1$, then $f$ is strictly increasing and hence it is invertible. This implies that $\gamma$ encloses a star-shaped domain around the origin, and thus $\rho$ is constant by Corollary \ref{cor:star_shaped} and Lemma \ref{lem:gamma_gometry}.

If $k=0$, then
\begin{align*}
    f(u) = \frac{cu^*}{2}\left(\frac{u}{u^*}-\frac{I(u)}{I^*}\right) +f(0)
\end{align*}
and 
\begin{align*}
     \frac{df}{du} = \frac{cu^*}{2}\left(\frac{1}{u^*}-\frac{1}{I^*\rho(u)^2}\right) = \frac{c}{2}\left(1-\frac{u^*}{I^*\rho(u)^2}\right).
\end{align*}
Plugging this into the constraint
\begin{align*}
     \left(\frac{\rho'(u)}{\rho(u)}\right)^2 +\left( f'(u)\right)^2 &= 1
\end{align*}
gives us

\begin{align*}
     \frac{\rho'(u)^2}{\rho(u)^2}  +\frac{c^2}{4}\left(1 - \frac{u^*}{I^*}\rho^{-2}\right)^2 &= 1\\
\end{align*}
which can also be written as
\begin{align*}
\rho^{4} \left(\frac{d}{d u} \left(\rho^{-2}\right)\right)^2+c^2\left(1 - \frac{u^*}{I^*}\rho^{-2}\right)^2 &= 4.    
\end{align*}
Let $r = \rho^{-2}$ and $\underline{r} = \frac{I^*}{u^*}$. Then, in terms of $r$, the constraint reads as
\begin{align*}
          r^{-2} \left(\frac{dr}{du}\right)^2+c^2\left(1 - \frac{r}{\underline{r}}\right)^2 &= 4.
\end{align*}
Rearranging for $\left(\frac{dr}{du}\right)^2$ gives us
\begin{align*}
     \left(\frac{dr}{du}\right)^2 &= 4r^2 - c^2r^2\left(1 - \frac{r}{\underline{r}}\right)^2.
\end{align*}
Divide the equation by $\underline{r}^2$ and let $R = \frac{r}{\underline{r}}$. Then we get
\begin{align*}
     (R')^2&=4R^2 - c^2R^2(R-1)^2\\
  &=c^2R^2 \left(\frac{\abs{c}-2}{\abs{c}}-R\right)\left( R - \frac{\abs{c}+2}{\abs{c}}\right).
\end{align*}
Since $c = \frac{2\pi}{u^*}>0$ by (\ref{eqn:c_u_star}), we have
\begin{align*}
     (R')^2&=c^2R^2 \left(\frac{c-2}{c}-R\right)\left( R - \frac{c+2}{c}\right).
\end{align*}
From this equation, it is clear that $R$ can have at most two critical values. Also, the winding number of $\gamma$ around the origin is $k=0$, so $\rho$ cannot be constant, which in turn implies that $R$ is not constant. Since $R$ is a non-constant periodic function, its global minimum and maximum values are distinct and are both attained. Hence, these are the only possible critical values, and we can write
\begin{align*}
     \left(R'\right)^2&=c^2R^2 \left(R_{min}-R\right)\left( R - R_{max}\right)
\end{align*}
where $R_{min} = \frac{c-2}{c}$ and $R_{max} = \frac{c+2}{c}$. Since $R$ is a strictly positive periodic function, we must have $R_{min}>0$, which in turn implies that $c >2$. 
We can also conclude that every critical point of $R$ yields either a global minimum or a global maximum. Since $R$ is just $\rho^{-2}$ normalized, we also get that every critical point of $\rho$ is either a global minimum or a global maximum point. Since $\gamma$ does not enclose the origin, by Lemma \ref{lemma:crit_points}, $\rho$ has exactly two critical points in $[0,u^*)$. Therefore, $R$ has exactly two critical points as well corresponding to its global extrema.  Without the loss of generality, let us assume that $R(0) = R_{max}$ and let $u_1$ denote the minimum point of $R$, i.e. $R(u_1) = R_{min}$. Then on $(0,u_1)$, $R' <0$ and on $(u_1,u^*)$ $R'>0$. So for $u \in (0,u_1)$,

\begin{align*}
    c = - \frac{R'}{R\sqrt{\frac{4}{c^2}-(R-1)^2}}
\end{align*}
and
\begin{align*}
    cu  &= -\int_0^u \frac{R'(z)}{R(z)\sqrt{\frac{4}{c^2}-(R(z)-1)^2}} dz\\
                    &= -\frac{c}{2}\int_{\frac{c}{2}\left(R(0)-1\right)}^{\frac{c}{2}\left(R(u)-1\right)} \frac{d\delta}{(\delta +\frac{c}{2})\sqrt{1-\delta^2}}\\
                    &= - \frac{c}{2}\frac{2}{\sqrt{c^2-4}}\arctan\frac{c\delta+2}{\sqrt{c^2-4}\sqrt{1-\delta^2}}\biggr|_{\frac{c}{2}\left(R(0)-1\right)}^{\frac{c}{2}\left(R(u)-1\right)}\\
\end{align*}
where we used the change of variable $\delta = \frac{c}{2}\left(R(u)-1\right)$.
Since 
$$
\frac{c}{2}\left(R(0)-1\right) = \frac{c}{2}\left(R_{max}-1\right)= \frac{c}{2}\left(\frac{c+2}{c}-1\right) = 1,
$$ 
we get
\begin{align*}
    u&=\frac{1}{\sqrt{c^2-4}}\left(\frac{\pi}{2} - \arctan\frac{\frac{c^2}{2}\left(R(u)-1\right)+2}{\sqrt{c^2-4}\sqrt{1-\frac{c^2}{4}\left(R(u)-1\right)^2}}\right).
\end{align*}
Also, since 
$$
\frac{c}{2}\left(R(u_1)-1\right) =\frac{c}{2}\left(R_{min}-1\right) =\frac{c}{2}\left(\frac{c-2}{c}-1\right) = -1
$$ 
and $c>2$, we get
\begin{align*}
    u_1&=\frac{1}{\sqrt{c^2-4}}\left[\frac{\pi}{2} - \left(-\frac{\pi}{2}\right)\right] = \frac{\pi}{\sqrt{c^2-4}}.
\end{align*}
Similarly, for $u \in(u_1,u^*)$,
\begin{align*}
    c =  \frac{R'}{R\sqrt{\frac{4}{c^2}-(R-1)^2}}
\end{align*}
and
\begin{align*}
    cu -cu_1  &= \int_{u_1}^u \frac{R'(z)}{R(z)\sqrt{\frac{4}{c^2}-(R(z)-1)^2}} dz\\
                    &=  \frac{c}{2}\frac{2}{\sqrt{c^2-4}}\arctan\frac{c\delta+2}{\sqrt{c^2-4}\sqrt{1-\delta^2}}\biggr|_{\frac{c}{2}\left(R(u_1)-1\right)}^{\frac{c}{2}\left(R(u)-1\right)}.
\end{align*}
Again, since $\frac{c}{2}\left(R(u_1)-1\right) = -1$, and $c>2$, we get
\begin{align*}
    u -u_1&=\frac{1}{\sqrt{c^2-4}}\left[\arctan\frac{\frac{c^2}{2}\left(R(u)-1\right)+2}{\sqrt{c^2-4}\sqrt{1-\frac{c^2}{4}\left(R(u)-1\right)^2}}- \left(-\frac{\pi}{2}\right)\right]\\
    &=\frac{1}{\sqrt{c^2-4}}\left[\arctan\frac{\frac{c^2}{2}\left(R(u)-1\right)+2}{\sqrt{c^2-4}\sqrt{1-\frac{c^2}{4}\left(R(u)-1\right)^2}} + \frac{\pi}{2}\right].
\end{align*}
Since 
$$
\frac{c}{2}\left(R(u^*)-1\right) = \frac{c}{2}\left(R_{max}-1\right)= \frac{c}{2}\left(\frac{c+2}{c}-1\right) = 1,
$$ 
we get
\begin{align*}
    u^* - u_1&=\frac{1}{\sqrt{c^2-4}}\left[\frac{\pi}{2} + \frac{\pi}{2}\right] = \frac{\pi}{\sqrt{c^2-4}}.
\end{align*}
Therefore,
\begin{align*}
     u^* &= \frac{2 \pi}{\sqrt{c^2-4}} > \frac{2\pi}{c}
\end{align*}
which contradicts (\ref{eqn:c_u_star}).

\appendix
\section{Construction of $L_\gamma$ by means of symplectic reduction} \label{appendix:construction}
 The twist tori $L_\gamma$ can be constructed by means of symplectic reduction, see \citep{brendel2021introduction}. This construction involves the identification of a certain reduced symplectic space by $(\C^*,\omega_\C)$ where $\omega_\C = \frac{i}{2}dw\wedge d\bar{w}$ is the standard symplectic form on $\C^*$. The goal of this Appendix is to discuss this identification in detail. 
 
Consider $\C^2$ equipped with the standard symplectic form $\omega_{\C^2} = \frac{i}{2}\left(dz_1 \wedge d \bar{z}_1 + dz_2 \wedge d \bar{z}_2\right)$. Let $h(z_1,z_2):\C^2 \to \R$ and $l:\C^2 \to \C$ be given by
\[
    h(z_1,z_2) = \frac{1}{2}\left(\abs{z_1}^2 - \abs{z_1}^2\right)
\]
and
\[
l(z_1,z_2) = z_1z_2.
\]

Then the flow of the vector field $X_h=J \nabla h$ generates an $S
^1$-action on $\C^2$. This action $A:S^1 \times \C^2 \to \C^2$ is given by 
\[
    A[t,(z_1,z_2)] = \left(z_1 e^{it},z_2 e^{-it}\right).
\]
The function $h$ is invariant under the action $A$ and thus it can be restricted to its level sets. In particular, $A$ can be restricted to $Z_0 := h^{-1}(0)$. Let $M_0^*$ be the quotient of $Z_0^* :=Z_0 \setminus \{0\}$ by this action and let $\pi_0 :Z_0^* \to M_0^* $ be the quotient map. The differential form $\omega_{\C^2}|_{Z_0}$ is invariant under $A$ so there exists a closed $2$-form $\omega_0$ on $M^*_0$ that satisfies $\omega_{\C^2}|_{Z_0^*} =\pi_0^*\omega_0$. Moreover, $\omega_0$ is non-degenerate so it is a symplectic form. The symplectic manifold $(M_0^*,\omega_0)$ is an example of a symplectic reduced space. For a more detailed discussion see \citep{brendel2021introduction}.

The function $l$ is also invariant under $A$ so it descends to $M_0^*$ and there is a well-defined smooth function $\Tilde{l}:M_0^* \to \C^*$ that satisfies
\[
    l|_{Z_0^*} = \Tilde{l} \circ \pi_0.
\]
By the chain rule, we also have
\[ 
    d\left(l|_{Z_0^*}\right) = d\Tilde{l} \circ d\pi_0.   
\]
 The quotient map $\pi_0$ is a submersion so $d\pi_0$ is surjective which tells us that $d\left(l|_{Z_0^*}\right)$ is surjective if and only if $d\Tilde{l}$ is surjective. Also, since $\dim_\R M_0^* =  \dim_\R \C^* =2$, the linear map $d\tilde{l}$ is bijective if and only if it is surjective. Therefore, $d\tilde{l}$ is bijective if and only if $d\left(l|_{Z_0^*}\right)$ is surjective. The level set $Z_0$ is just the cone over the torus $S^1(1) \times S^1(1) \subset \C^2$ so $Z_0^*$ can be parametrized by $G:(0,\infty) \times [0,2\pi) \times [0,2\pi) \to \C^2$ that is given by
 \begin{align*}
     G(r,\theta,\eta) = r\left(e^{i\theta},e^{i\eta}\right).
 \end{align*}
 In these coordinates, $l|_{Z_0^*}$ is given by
 \begin{align*}
     (r,\theta,\eta) \mapsto r^2e^{i(\theta+\eta)}
 \end{align*}
 so
 \begin{align*}
      d\left(l|_{Z_0^*}\right) = 2re^{i(\theta+\eta)}dr + ir^2e^{i(\theta+\eta)}\left(d\theta + d\eta\right).
 \end{align*}
Fix a point $(r_0,\theta_0,\eta_0) \in Z_0^*$ and consider a tangent vector at $(r_0,\theta_0,\eta_0)$ of the form $X=A\partial_r + B\partial_\theta$ where $A$ and $B$ are real constants. Then
\begin{align*}
     d\left(l|_{Z^*_0}\right) (X) &= 2r_0e^{i(\theta_0+\eta_0)}A + i\left(r_0\right)^2e^{i(\theta_0+\eta_0)}B\\
     &=r_0e^{i(\theta_0+\eta_0)}\left(2A +ir_0B\right).
\end{align*}
For $z \in \C$, let $\frac{z}{2r_0e^{i(\theta_0+\eta_0)}} = x+iy$. Then the equation $d\left(l|_{Z^*_0}\right) (X) = z$ is equivalent to
\begin{align*}
    2A +ir_0B = x+iy
\end{align*}
which is clearly solved by $A = \frac{x}{2}$ and $B=\frac{y}{r_0}$. Therefore, we can conclude that $d\left(l|_{Z^*_0}\right)$ is surjective at all points of $Z_0^*$. By our previous discussion, $d\tilde{l}$ must be bijective and hence $\tilde{l}$ is a local diffeomorphism. 

Next, we show that $\Tilde{l}$ is a bijection. It is clearly surjective as for any $w = \rho e^{i\varphi} \in \C^*$ we have 
\[
    w = l\left(\sqrt{\rho}e^{i\frac{\varphi}{2}},\sqrt{\rho}e^{i\frac{\varphi}{2}}\right) = \tilde{l}\left(\pi_0\left(\sqrt{\rho}e^{i\frac{\varphi}{2}},\sqrt{\rho}e^{i\frac{\varphi}{2}}\right)\right),
\]
i.e. 
\begin{align*}
    w= \tilde{l}(z)
\end{align*}
for $z =\pi_0\left(\sqrt{\rho}e^{i\frac{\varphi}{2}},\sqrt{\rho}e^{i\frac{\varphi}{2}}\right) \in M_0^*$.
To prove injectivity, let us assume that $\tilde{l}(\pi_0(z_1,z_2)) = \tilde{l}(\pi_0(v_1,v_2))$ for some $(z_1,z_2),(v_1,v_2) \in Z_0^*$. We need to show that $\pi_0(z_1,z_2) = \pi_0(v_1,v_2)$ or equivalently that
\[
    (v_1,v_2) =\left(z_1 e^{it},z_2 e^{-it}\right)
\]
for some $t \in S^1$. Since, $\tilde{l}(\pi_0(z_1,z_2)) = \tilde{l}(\pi_0(v_1,v_2))$,
\[
    v_1v_2 = z_1z_2.
\]
Also, since $(v_1,v_2),(z_1,z_2) \in Z_0$, we have $\abs{v_1} = \abs{v_2}$ and $\abs{z_1} = \abs{z_2}$. Therefore, we can write $v_1 = r_ve^{i\theta_v},v_2 = r_ve^{i\eta_v},z_1 = r_ze^{i\theta_z}$ and $z_2 = r_ze^{i\eta_z}$ and have
\[
    r_v^2e^{i(\theta_v + \eta_v)} = r_z^2e^{i(\theta_z + \eta_z)}.
\]
Thus $r_v=r_z$ and 
\begin{align*}
      \theta_v + \eta_v &\equiv \theta_z + \eta_z \mod 2\pi
\end{align*}
or equivalently,
\begin{align*}
      \theta_v -\theta_z   &\equiv  \eta_z -\eta_v \mod 2\pi.
\end{align*}
So $e^{-i( \theta_v -\theta_z) } = e^{-i( \eta_z -\eta_v)}$ and for $t = \theta_v -\theta_z \pmod{2 \pi}$,
\begin{align*}
    \left(z_1 e^{it},z_2 e^{-it}\right) &= \left(r_ze^{i\theta_z}e^{i( \theta_v -\theta_z)},r_ze^{i\eta_z} e^{-i( \eta_z - \eta_v)}\right) \\
                                     =&\left(r_ve^{i\theta_v},r_ve^{i\eta_v}\right)\\
                                     =&(v_1,v_2).
\end{align*}
Therefore, $\Tilde{l}$ is a bijection, and we can conclude that it is a diffeomorphism.

$Z_0$ is given by the equation
\begin{align*}
    \abs{z_1}^2 = \abs{z_2}^2
\end{align*}
or equivalently,
\begin{align*}
    z_1\Bar{z}_1 = z_2 \Bar{z}_2.
\end{align*}
Differentiating this equation gives us
\begin{align*}
    \bar{z}_1dz_1 + z_1d\bar{z}_1 &= \bar{z}_2dz_2 + z_2d\bar{z}_2.
\end{align*}
Multiplying both sides by $z_1$, we get
\begin{align*}
    \abs{z_1}^2 dz_1 + z_1^2d\bar{z}_1&=z_1\bar{z}_2dz_2 + z_1z_2d\bar{z}_2
\end{align*}
so along $Z_0^*$,
\begin{align*}
     \abs{z_1}^2 dz_1 \wedge d\bar{z}_1&=z_1\bar{z}_2dz_2 \wedge d\bar{z}_1 + z_1z_2d\bar{z}_2 \wedge d\bar{z}_1.
\end{align*}
Similarly,
\begin{align*}
     \abs{z_2}^2 dz_2 \wedge d\bar{z}_2&=z_2\bar{z}_1dz_1 \wedge d\bar{z}_2 + z_1z_2d\bar{z}_1 \wedge d\bar{z}_2.
\end{align*}
Therefore,
\begin{align*}
     \omega_{\C^2}|_{Z_0^*}&=\frac{i}{2}\left(dz_1 \wedge d\bar{z}_1+dz_2 \wedge d\bar{z}_2\right)\\
     &=\frac{i}{2\abs{z_1}^2}\left(z_2\bar{z}_1dz_1 \wedge d\bar{z}_2 + z_1\bar{z}_2dz_2 \wedge d\bar{z}_1\right)
\end{align*}
and the pull-back of the symplectic form $\omega := \frac{1}{2\abs{w}}\omega_\C  = \frac{i}{4\abs{w}}dw \wedge d\bar{w}$ by $l|_{Z_0^*}$ is given by
\begin{align*}
    \left(l|_{Z_0^*}\right)^*\omega&=\frac{i}{4\abs{l}}dl \wedge d\bar{l}\\
    &= \frac{i}{4\abs{z_1z_2}}\left(z_1 dz_2 + z_2 dz_1\right) \wedge \left(\bar{z}_1 d\bar{z}_2 + \bar{z}_2 d\bar{z}_1\right)\\
    &= \frac{i}{4\abs{z_1} \abs{z_2}} \left(\abs{z_2}^2 dz_1 \wedge d\bar{z}_1 + \abs{z_1}^2 dz_2 \wedge d\bar{z}_2 + z_2\bar{z}_1dz_1 \wedge d\bar{z}_2 + z_1\bar{z}_2dz_2 \wedge d\bar{z}_1\right)\\
    &= \frac{i}{4}  \left( dz_1 \wedge d\bar{z}_1 + dz_2 \wedge d\bar{z}_2\right) +   \frac{i}{4\abs{z_1}^2}\left(z_2\bar{z}_1dz_1 \wedge d\bar{z}_2 + z_1\bar{z}_2dz_2 \wedge d\bar{z}_1\right) \\
    &=\frac{1}{2}\omega_{\C^2}|_{Z_0^*} + \frac{1}{2}\omega_{\C^2}|_{Z_0^*}\\
    &=\omega_{\C^2}|_{Z_0^*}
\end{align*}
Therefore, $\Tilde{l}$ is a symplectomorphism between $(M_0^*,\omega_0)$ and $(\C^*,\omega)$. In polar coordinates $w = \rho e^{i\varphi}$, the symplectic form $\omega$ is given by
\[
    \omega = \frac{1}{2} d\rho \wedge d\varphi.
\]
Let $\psi :\C^* \to \C^*$ be given by
\[
    \psi(w) = \abs{w}w
\]
or in polar coordinates by
\[
    \psi(\rho,\varphi) = (\rho^2,\varphi).
\]
Then $\psi$ is a diffeomorphism and the pull-back of $\omega$ by $\psi$ is 
\begin{align*}
    \psi^*\omega &= \frac{1}{2}d\left(\rho^2\right)\wedge d\varphi\\
    &=\rho d\rho\wedge d\varphi\\
    &=\omega_\C.
\end{align*}
We conclude that $(M_0^*,\omega_0)$ and $(\C^*,\omega_\C)$ are symplectomorphic via $\tilde{\psi} : = \psi^{-1} \circ \tilde{l}$.

Let $\gamma$ be a simple closed curve in $\C^*$ and define
\[
    L_{\gamma} : = \left\{\frac{1}{\sqrt{2}}\left(\gamma e^{i\alpha}, \gamma e^{-i\alpha}\right) \bigg| \hspace{5pt} \alpha \in [0,2\pi) \right\}.
 \]
We want to realize $L_\gamma$ as a lift of a curve $\tilde{\gamma} \subset (\C^*,\omega_\C)$, i.e. we are looking for a $\tilde{\gamma}$ for which
\begin{align*}
    L_\gamma = \pi_0^{-1}\left(\tilde{\psi}^{-1}\left(\tilde{\gamma}\right)\right).
\end{align*}
Such a $\tilde{\gamma}$ satisfies 
\begin{align*}
    \tilde{\gamma} &=\tilde{\psi}\left(\pi_0\left(L_\gamma\right)\right)\\  
            &=\psi^{-1}\left(\tilde{l}\left(\pi_0\left(L_\gamma\right)\right)\right)\\
            &=\psi^{-1}\left(l\left(L_\gamma\right)\right)\\
            &=\psi^{-1}\left(\frac{\gamma^2}{2}\right)\\
            &=\frac{\gamma^2}{\sqrt{2}\abs{\gamma}}.
 \end{align*}
Therefore, if $\gamma$  is parametrized by $\gamma(\beta)= \rho(\beta) e^{if(\beta)}$ where $\beta \in [0,2\pi)$, then $\tilde{\gamma}$ is given by 
\[
\tilde{\gamma}(\beta)= \frac{1}{\sqrt{2}}\rho(\beta) e^{i2f(\beta)}.
\]

When $\gamma$ is contained in the open upper half-plane, then it is possible to choose an identification so that $\tilde{\gamma} = \gamma$. Let $H:= \left\{w=\rho e^{i\varphi} \in \C: \rho>0, \varphi \in (0,\pi)\right\}$ denote the open upper half-plane and let $\C^*_{-} := \C^* \setminus \left\{w=\rho e^{i\varphi} \in \C:  \varphi = 0 \right\}$ denote the open set we get after removing the positive $x$-axis from $\C^*$. Let $\phi: H \to \C^*_{-}$ be defined by
\begin{align*}
    \phi(w) = \frac{1}{2}w^2.
\end{align*}
Then, in polar coordianates, $\phi$ is given by
\begin{align*}
    \phi(\rho,\varphi) = \left( \frac{1}{2} \rho^2,2\varphi \right)
\end{align*}
which is a diffeomorphism. Also,
\begin{align*}
    \phi^*\omega&= \frac{1}{2}d\left(\frac{1}{2} \rho^2\right) \wedge d(2 \varphi)\\
            &= \rho d\rho \wedge d\varphi\\
            &= \omega_\C.
\end{align*}
Therefore, $\phi$ is s symplectomorphism between $(H,\omega_\C)$ and $(\C^*_{-}, \omega)$. Suppose that 
$\gamma \subset H$. Then $l(L_\gamma) = \frac{1}{2}\gamma^2 \subset \C^*_{-}$ and $\phi^{-1}(l(L_\gamma)) = \gamma$. Let $(M_0^*)_{-} : = M^*_0\setminus \tilde{l}^{-1}\left(\left\{w=\rho e^{i\varphi} \in  \C:  \varphi = 0 \right\}\right)$. Then using the identification  $\tilde{\phi}:=\phi^{-1} \circ \tilde{l}: \left((M_0^*)_{-},\omega_0\right) \to \left(H,\omega_\C\right)$, we can write
\begin{align*}
    L_\gamma = \pi_0^{-1}\left(\tilde{\phi}^{-1}\left(\gamma\right)\right).
\end{align*}
\section{Embeddedness of $L_\gamma$}\label{appendix:embeddedness}
Let $\gamma \subset \C^*$ be a simple closed curve parametrized by $\gamma(\beta) = \rho(\beta) e^{if(\beta)}$. Consider the immersion $F:S^1 \times S^1 \to \C^2$ given by
\begin{align*}
    F(\alpha,\beta) = \frac{1}{\sqrt{2}}\left(\gamma e^{i\alpha},\gamma e^{-i\alpha}\right) = \frac{\rho(\beta)}{\sqrt{2}} \left(e^{i\left(f(\beta) +\alpha\right)},e^{i\left(f(\beta) - \alpha\right)}\right). 
\end{align*}
Then
\begin{align*}
    \partial_\alpha F &= i \frac{\rho}{\sqrt{2}} \left(e^{i\left(f +\alpha\right)},-e^{i\left(f - \alpha\right)}\right)\\
        \partial_\beta F &=\frac{\dot{\rho} + i \dot{f}\rho}{\sqrt{2}}\left(e^{i\left(f +\alpha\right)},e^{i\left(f - \alpha\right)}\right) = \frac{\dot{\gamma}}{\sqrt{2}}\left(e^{i\alpha},e^{-i\alpha}\right).
\end{align*}
Assuming that $\gamma$ is given a regular parametrization, $\partial_\beta F$ never vanishes and $\partial_\alpha F$ vanishes if and only if $\gamma$ passes through the origin. Therefore, when $\gamma$ does not pass through the origin, $L_\gamma := F(S^1 \times S^1)$ is an immersed torus.

Let us write the tangent vectors $\partial_\alpha F$ and $\partial_\beta F$ as
\begin{align*}
   \partial_\alpha F& = (iz_1,-iz_2)\\
   \partial_\beta F &= (\tau z_1,\tau z_2)
\end{align*}
where $\tau = \frac{\dot{\gamma}}{\gamma}$. Since $\L_\gamma \subset Z_0 = \left\{(z_1,z_2) \in  \C^2 :H(z_1,z_2) = \frac{1}{2}\left(\abs{z_1}^2 - \abs{z_2}^2\right) = 0\right\}$,
\begin{align*}
    \langle \partial_\beta F,J\partial_\alpha F\rangle & = \text{Re} \left[(\tau z_1,\tau z_2) \cdot \overline{(iz_1,-iz_2)}\right]\\
                                                        &= \text{Re} \left[-i\tau\left(\abs{z_1}^2 - \abs{z_2}^2\right)\right]\\
                                                        &=0
\end{align*}
and
\begin{align*}
    \langle \partial_\alpha F,J\partial_\beta F\rangle & =-  \langle \partial_\beta F,J\partial_\alpha F\rangle\\
                                                        &=0
\end{align*}
so $L_\gamma$ is a Lagrangian immersion.

Next, we examine the possible self-intersection points of $L_\gamma$. Suppose, that $F(\alpha_1,\beta_1) = F(\alpha_2,\beta_2)$ for some $(\alpha_1,\beta_1), (\alpha_2,\beta_2) \in S^1 \times S^1$. Then 
\begin{align*}
    \rho(\beta_1) = \abs{F(\alpha_1,\beta_1)} = \abs{F(\alpha_2,\beta_2)} = \rho(\beta_2)
\end{align*}
so we get the simplification
\begin{align*}
    F(\alpha_1,\beta_1) &= F(\alpha_2,\beta_2)\\
   \frac{ \rho(\beta_1)}{\sqrt{2}} \left(e^{i\left(f(\beta_1) +\alpha_1\right)},e^{i\left(f(\beta_1) - \alpha_1\right)}\right) &=  \frac{\rho(\beta_2)}{\sqrt{2}} \left(e^{i\left(f(\beta_2) +\alpha_2\right)},e^{i\left(f(\beta_2) - \alpha_2\right)}\right)\\
    \left(e^{i\left(f(\beta_1) +\alpha_1\right)},e^{i\left(f(\beta_1) - \alpha_1\right)}\right) &= \left(e^{i\left(f(\beta_2) +\alpha_2\right)},e^{i\left(f(\beta_2) - \alpha_2\right)}\right).
\end{align*}
Identifying $S^1$ with $\R / [0,2\pi]$, we see that this equation holds if and only if
\begin{align*}
    f(\beta_1) +\alpha_1 &\equiv f(\beta_2) +\alpha_2 \mod 2\pi\\
    f(\beta_1) - \alpha_1 &\equiv f(\beta_2) - \alpha_2 \mod 2\pi.
\end{align*}
Adding these two equations tells us that
\begin{align*}
    2\left(f(\beta_1) - f(\beta_2) \right)  &\equiv 0 \mod 2\pi
\end{align*}
or equivalently, that 
\begin{align*}
    f(\beta_1) = f(\beta_2) + n\pi
\end{align*}
for some $n \in \Z$. There are two solutions to this equation in $S^1$, namely, $f(\beta_1) = f(\beta_2)  \pmod{ 2 \pi}$ and $f(\beta_1) = f(\beta_2) +\pi \pmod{2 \pi}$. In the first case, we have $\gamma(\beta_1) = \gamma(\beta_2)$ which is impossible, since $\gamma$ is a simple curve. Therefore, the only possible solution is given by 
\begin{align}\label{eqn:appendix_f}
    f(\beta_1) = f(\beta_2) +\pi \mod2 \pi
\end{align}
which also tells us that
\begin{align}\label{eqn:appendix_alpha}
    \alpha_1 \equiv \alpha_2 +\pi \mod{2 \pi}.
\end{align}
Since $\rho(\beta_1) = \rho(\beta_2)$ and $ f(\beta_1) = f(\beta_2) +\pi \pmod{2 \pi}$, we see that $\gamma(\beta_1)$ and $\gamma(\beta_2)$ are exactly the two intersection points of $C_{\rho(\beta_1)}$, the circle of radius $\rho(\beta_1)$ centred at the origin, and the line $\{w=\rho e^{i\varphi} \in \C : r \in (-\infty,\infty), \varphi = \beta_1\}$. Since $\gamma$ is simple, it passes through both $\gamma(\beta_1)$ and $\gamma(\beta_2)$ exactly once. Hence $\beta_1$ determines $\beta_2$ uniquely and vice versa. By equation (\ref{eqn:appendix_alpha}), $\alpha_1$ also determines $\alpha_2$ uniquely and vice versa. Therefore, every self-intersection point of $L_\gamma$ is a double point.

Let $\underline{z} \in L_\gamma$ be a double point such that $\underline{z} = F(\alpha_1,\beta_1) = F(\alpha_2,\beta_2)$. Then, by (\ref{eqn:appendix_f}) and (\ref{eqn:appendix_alpha}) the tangent basis vectors at $\underline{z}$ induced by the parametrization are given by
\begin{align*}
    \partial_\alpha  F(\alpha_1,\beta_1) &= i \rho(\beta_1) \left(e^{i\left(f(\beta_1) +\alpha_1\right)},-e^{i\left(f(\beta_1) - \alpha_1\right)}\right) \\
    &= i \rho(\beta_2) \left(e^{i\left(f(\beta_2) + \pi +\alpha_2 + \pi\right)},-e^{i\left(f(\beta_2) +\pi - \alpha_2 - \pi\right)}\right)\\
    &= i \rho(\beta_2) \left(e^{i\left(f(\beta_2) +\alpha_2\right)},-e^{i\left(f(\beta_2) - \alpha_2\right)}\right)\\
    &= \partial_\alpha  F(\alpha_2,\beta_2), \\
     \partial_\beta  F(\alpha_1,\beta_1) &= \dot{\gamma}(\beta_1)\left(e^{i\alpha_1},e^{-i\alpha_1}\right)
\end{align*}
and
\begin{align*}
     \partial_\beta  F(\alpha_2,\beta_2) &= \dot{\gamma}(\beta_2)\left(e^{i\alpha_2},e^{-i\alpha_2}\right)\\
                                        &= \dot{\gamma}(\beta_2)\left(e^{i\alpha_1 +i\pi},e^{-i\alpha_1 -i\pi}\right)\\
                                        &=-\dot{\gamma}(\beta_2)\left(e^{i\alpha_1},e^{-i\alpha_1}\right).
\end{align*}
Therefore, $L_\gamma$ "touches" itself at $\underline{z}$ if and only if $\partial_\alpha  F(\alpha_1,\beta_1), \partial_\beta  F(\alpha_1,\beta_1)$ and $\partial_\beta  F(\alpha_2,\beta_2)$ span a $2$-plane. In Section \ref{section:DiffGeo}, we saw that $\partial_\alpha  F$ is always perpendicular to $\partial_\beta  F$ so $\partial_\beta  F(\alpha_1,\beta_1)$ is perpendicular to $\partial_\alpha  F(\alpha_1,\beta_1)$ and $\partial_\beta  F(\alpha_2,\beta_2)$ is perpendicular to $\partial_\alpha  F(\alpha_2,\beta_2) = \partial_\alpha  F(\alpha_1,\beta_1)$. Therefore, these three vectors span a $2$-plane if and only if $\partial_\beta  F(\alpha_1,\beta_1)$ and $\partial_\beta  F(\alpha_2,\beta_2)$ are parallel or equivalently if  $\dot{\gamma}(\beta_1)$ and $\dot{\gamma}(\beta_2)$ are parallel. Otherwise, $L_\gamma$ "crosses" itself. Since $ \partial_\alpha  F(\alpha_1,\beta_1)= \partial_\alpha  F(\alpha_2,\beta_2) $, the intersection is not transversal even in this case.  

\bibliographystyle{amsplain}

\bibliography{biblio}

\end{document}